\documentclass[a4paper, twocolumn]{revtex4}
\usepackage{amsmath, amssymb, txfonts}
\usepackage{graphicx}
\usepackage{fancyhdr}
\begin{document}
\title{Numerical analysis with the phase field equations for the Stefan problems}
\author{Jun-ichi Koga}
%\affiliation{Division of materials science, Graduate School of Science and Engineering, Saitama University, Sakura-ku, Saitama, 338-8570 Japan}
\begin{abstract}
Phase field equations describe the novel approach to the Stefan problems. We calculate these equations numerically performed in two-dimensions. We take full advantage of the phase field parameter $\varphi$ to track the interface on which complicated movements or changes are needed. More precisely, we analyze the dissolution and/or precipitation without the unknown velocity functions like the Navier-Stokes equations numerically.\\\\
PACs numbers: v79B 68. 70. w+.
\end{abstract}
\maketitle
\section{introduction}
In these thirty years, it seems that phase field equations describe the novel approach to the Stefan problems. For example, it has been proved that the phase field equations converge to the equations of the Stefan problem as a singular limit. Caginalp and a co-researcher [1-4] have studied the reduction of these equations in the singular limit of Stefan problems and Hele-Shaw problem. The solution of the Stefan problem is obtained from the phase-field equations which consist the heat equation and the Ginzburg-Landau equation by Soner [10]. Hence, by the phase-field equations, we determine the sharp profile of the phase, the temperature distribution and the velocity of the moving interface between a solid and a liquid.\\
\quad In this paper, using the Landau-Ginzburg theory of phase transition, physical properties show the phase field equations we consider, which is reduced to by the free energy
\begin{equation*}
F(\varphi)=\int_\Omega dx\left[\frac{\epsilon^2}{2}\left|\nabla\varphi\right|^2+W(\varphi)\right].
\end{equation*}
\quad First we use the phase field equations of ``Caginalp type'' due to [2],
\begin{equation}
\begin{array}{l}
\displaystyle u_t+\frac{\ell}{2}\varphi_t =\Delta u\\\\
\displaystyle\varphi_t=\Delta\varphi +\frac{1}{2}\left(\varphi -\varphi^3\right)+2u
\end{array}
\end{equation}
directly and numerically. Here $W(\varphi)=:\frac{1}{2}\left(\varphi -\varphi^3\right)$ is a derivative of a symmetric double-well potential with minima at $\pm 1$, $\ell$ a latent heat and $\varphi$ $+1$ liquid (e.g. water), $-1$ solid (e.g. ice).\\
\quad Second we investigate the following equations
\begin{equation}
\frac{\partial\varphi}{\partial t}=M_\varphi\left[\Delta\varphi+4g(\varphi)\varphi\left(1-\varphi\right)\left(\varphi-0.5+\beta\right)\right]
\end{equation}
in two space dimensions, where $M_\varphi$ is a mobility of $\varphi$ and these equations are called the modified Allen-Cahn equations.\\
\quad Third we consider the following equations due to [6],
\begin{equation}
\begin{array}{l}
\displaystyle \tau\frac{\partial\varphi}{\partial t}=W^2\frac{\partial^2\varphi}{\partial x^2}+\left[\varphi -\lambda u\left(1-\varphi^2\right)\right]\left[1-\varphi^2\right],\\\\
\displaystyle\frac{\partial u}{\partial t}=D\frac{\partial^2 u}{\partial x^2}+\frac{1}{2}\frac{\partial\varphi}{\partial t}
\end{array}
\end{equation}
Then applying the moving frame to these equations yields
\begin{equation}
\begin{array}{l}
\displaystyle\tau V\frac{\partial\varphi}{\partial x}+W^2\frac{\partial^2\varphi}{\partial x^2}+\left[\varphi -\lambda u\left(1-\varphi^2\right)\right]\left[1-\varphi^2\right]=0\\\\
\displaystyle V\frac{\partial u}{\partial x}+D\frac{\partial^2 u}{\partial x^2}-\frac{V}{2}\frac{\partial\varphi}{\partial x}=0
\end{array}
\end{equation}
\quad Finally, we also use the following relations due to [8],
\begin{equation}
\begin{array}{l}
\displaystyle\frac{\partial\varphi}{\partial t}=\frac{1}{Pe}\left[\nabla^2\varphi +\left(1-\varphi^2\right)\left(\varphi -\lambda c\right)-\kappa |\nabla\varphi|\right]\\\\
\displaystyle\frac{\partial c}{\partial t}=\nabla^2 c+\alpha\frac{\partial\varphi}{\partial t}+\left(\nabla^2\varphi-\frac{\partial\varphi}{\partial t}\right)\frac{\alpha\cdot\displaystyle\frac{\partial\varphi}{\partial t}}{Da|\nabla\varphi|}
\end{array}
\end{equation}
to simulate the symmetric circular domain of solid. Here $Pe$ is P\'eclet number, $\lambda$ a positive parameter, $\kappa$ a mean curvature, $\alpha$ also a positive parameter, $Da$ a Damk\"ohler number.\\
\quad The so-called one phase Stefan problem has no discontinuous-gradient model for the temperature distribution. However, in the two phase Stefan problem, this phenomena occurs. Thus one can obtain the sharp-interface numerically by investigating what control it.\\
\quad Our first purpose of this paper is to find out what determine the sharpness of the interface between two different phases (e.g. ice and water). The second is to investigate asymptotic behavior of the volume of the solidification and dissolution and/or precipitation.\\
\quad We calculate these equations numerically performed in two-dimensions. We take full advantage of the phase field parameter $\varphi$ to track the interface on which complicated movements or changes are needed. More precisely, we analyze the dissolution and/or precipitation without the unknown velocity functions like the Navier-Stokes equations numerically.
\section{numerical method}
In our survey, we use the finite differential method for the time evolution equations numerically.\\
\quad For Eqs. (1), we use the 5-point finite difference stencil
\begin{equation}
\Delta\varphi_{(i,j)}=\frac{\varphi_{(i+1,j)}+\varphi_{(i-1,j)}+\varphi_{(i,j+1)}+\varphi_{(i,j-1)}-4\varphi_{(i,j)}}{(\Delta x)^2}
\end{equation}
and we use the forward scheme in time
\begin{equation}
\frac{\partial\varphi}{\partial t}\cong\frac{\varphi^{(n+1)}_{(\ell, m)}-\varphi^{(n)}_{(\ell, m)}}{\Delta t},
\end{equation}
\quad For Eq. (2), we use the forward scheme in time
\begin{equation}
\frac{\partial\varphi}{\partial t}\cong\frac{\varphi^{(n+1)}_{(\ell, m)}-\varphi^{(n)}_{(\ell, m)}}{\Delta t},
\end{equation}
on the other hand, the second order center-ward scheme in space
\begin{equation}
\begin{array}{l}
\displaystyle\frac{\partial^2\varphi}{\partial x^2}\cong\frac{\varphi^{(n)}_{(\ell-1, m)}-2\varphi^{(n)}_{(\ell, m)}+\varphi^{(n)}_{(\ell+1,m)}}{(\Delta x)^2},\\\\
\displaystyle\frac{\partial^2\varphi}{\partial y^2}\cong\frac{\varphi^{(n)}_{(\ell, m-1)}-2\varphi^{(n)}_{(\ell, m)}+\varphi^{(n)}_{(\ell, m+1)}}{(\Delta y)^2}
\end{array}
\end{equation}
\quad For Eqs. (3), the gradient of $\varphi$ with respect to $x$, $y$, we use the central difference approximation
\begin{equation}
\begin{array}{l}
\displaystyle|\nabla\varphi|_{(i,j)}=\frac{1}{\Delta x}\sqrt{\frac{\left(\varphi_{(i+1,j)}-\varphi_{(i-1,j)}\right)^2}{4}+\frac{\left(\varphi_{(i,j+1)}-\varphi_{(i,j-1)}\right)^2}{4}},\\\\
\displaystyle|\nabla\varphi|_{(i,j)}=\frac{1}{\Delta y}\sqrt{\frac{\left(\varphi_{(i+1,j)}-\varphi_{(i-1,j)}\right)^2}{4}+\frac{\left(\varphi_{(i,j+1)}-\varphi_{(i,j-1)}\right)^2}{4}},
\end{array}
\end{equation}
respectively, and
\begin{equation}
\nabla^2\varphi_{(i,j)}=\frac{\varphi_{(i+1,j)}+\varphi_{(i-1,j)}+\varphi_{(i,j+1)}+\varphi_{(i,j-1)}-4\varphi_{(i,j)}}{(\Delta x)^2}.
\end{equation}
If we do not have any confusion, we abbreviate the ``$(n)$''.
\section{asymptotic behavior of solutions to Stefan problems}
This result is well-known. However we consider it to show an asymptotic behavior of the solutions for the Stefan problems.\\
\quad In general, the parabolic partial differential equations often have initial-boundary value problems in a certain sense and the exact solutions are difficult to obtain by the direct calculations. Therefore the fact that we investigate the asymptotic behavior of solutions of general case plays an important r\^ole in this sense. Thus we get
\begin{equation}
s(t)\sim\beta t^{1/2}
\end{equation}
where $\beta$ corresponds to the integral formula ``Gaussian distribution'' (see Friedman, 1964, 2008).
\section{results and discussions}
\begin{center}
TABLE I. Numerical volumes of the solidification and dissolution and/or precipitation for the time step. From the above, Eq. (1) v.s. Volume 1 and Eq. (2) v.s. Volume 2.
\end{center}
\begin{center}
\begin{tabular}{cccccccc}
\hline
nstep & 20 & 40 & 60 & 80 & 100 &\\
\hline
\hline
Volume 1 & 4.8833 & 4.8848 & 4.8871 & 4.8899 & 4.8931 & $(\times 10^4)$ & \\
Volume 2 & 3.5925 & 4.0135 & 4.4569 & 4.9166 & 5.3904 &  $(\times 10^4)$ &\\
\hline
\end{tabular}
\end{center}
\hspace{10truemm}\\
\quad We obtain the fact that the interface between two regions in the solidification and dissolution and/or precipitation are determined by the ``coupling parameters'' and coefficients, and sharp in our numerical simulation.\\
\quad The volume controlled by the solidification and dissolution and/or precipitation is the dimensionless versions. we report that the volume of them is proportional to the time step (nstep in TABLE I).
\begin{figure}[h]
\begin{center}
\includegraphics[width=10cm, height=13cm, clip]
{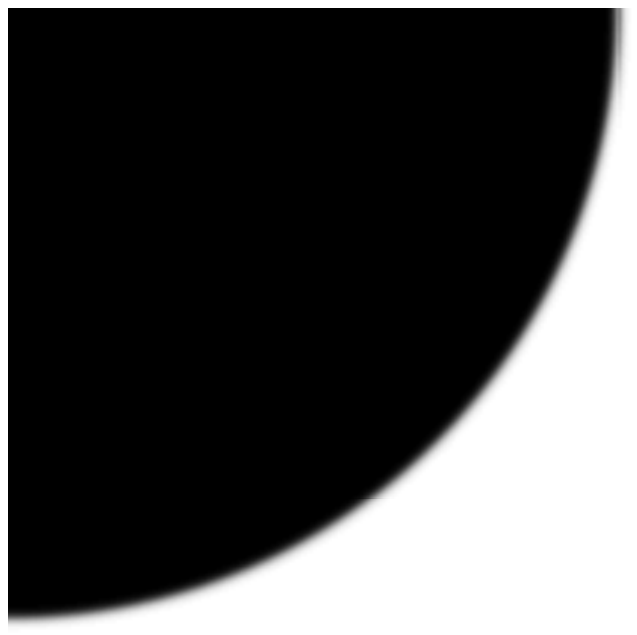}
\caption{Profile of $\phi$.}
\end{center}
\end{figure}

\end{document}